\documentclass[12pt]{amsart}
\usepackage{amsmath, amsthm}
\newtheorem{thm}{Theorem}[section]
\newtheorem{lemma}[thm]{Lemma} 

\theoremstyle{definition}
\newtheorem*{defn}{Definition}
\newenvironment{defncite}[1]{\begin{defn}[#1]}{\end{defn}}
\newtheorem{ex}[thm]{Example}

\theoremstyle{remark}

\newcounter{alphcount}
\renewcommand{\thealphcount}{(\alph{alphcount})}
\newenvironment{alph-list}{\begin{list}{\thealphcount}{\usecounter{alphcount}}}%
                        {\end{list}} 

\newcommand{\ie}{{\em i.e.}} 
\newcommand{\eg}{{\em e.g.}} 
\newcommand{\dfn}[1]{{\bf {#1}}} 
 
\newcommand{\abs}[1]{|#1|}

\newcommand{\disun}{\mathbin{\dot{\cup}}}       %from Mark Purtill
\newcommand{\condns}[2]{\substack{#1 \\ #2}}
\newcommand{\numsys}[1]{{\bf #1}}
\renewcommand{\colon}{{:}\;}
\newcommand{\maps}{\rightarrow}

\newcommand{\rhom}{\widetilde{H}}
\newcommand{\rhomi}[1]{\widetilde{H}_{#1}}

\newcommand{\nd}{\;}    %''and'' in series of conditions

\newcommand{\seq}[3]{{#1}_{#2},\dots,{#1}_{#3}}
\newcommand{\ltseq}[3]{{#1}_{#2} < \dots < {#1}_{#3}}
\newcommand{\ltsubseq}[4]{#1_{#2_#3} < \dots < #1_{#2_#4}}

\newcommand{\set}[3]{\{\seq{#1}{#2}{#3}\}}
\newcommand{\setseq}[3]{\{\ltseq{#1}{#2}{#3}\}}
\newcommand{\setsubseq}[4]{\{\ltsubseq{#1}{#2}{#3}{#4}\}}

\newcommand{\wedgesubseq}[4]{#1_{#2_#3} \wedge \dots \wedge #1_{#2_#4}} 

\newcommand{\kk}{\boldsymbol k}
\newcommand{\dimk}{\dim_{\kk}}
\newcommand{\kV}{\kk V}
\newcommand{\imt}[1]{\widetilde{#1}}%''IMage-Tilde''
\newcommand{\imf}[1]{\widetilde{f}_{#1}}%''IMage-Tilde'' -- f-variation
\newcommand{\img}[1]{\widetilde{\gamma}_{#1}}%...gamma variation
\newcommand{\ime}[1]{\widetilde{e}_{#1}}%...e variation
\newcommand{\collC}{\mathcal C}

\DeclareMathOperator{\spn}{span}
\DeclareMathOperator{\im}{im}
\DeclareMathOperator{\Gin}{Gin}

\newcommand{\betat}{\beta}
\newcommand{\rbeti}[1]{\beta_{#1}}

\begin{document} 

\title{Algebraic Shifting Increases Relative Homology}

\author{Art M. Duval}
\address{\hskip-\parindent
        University of Texas at El Paso\\
        Department of Mathematical Sciences\\   
        El Paso, TX 79968-0514}
\email{artduval@math.utep.edu}

\thanks{Research performed while a General Member
of the Mathematical Sciences Research Institute, Berkeley.  Research
at MSRI is supported in part by NSF grant DMS-9701755}

\begin{abstract}

We show that algebraically shifting a pair of simplicial complexes
weakly increases their relative homology Betti numbers in every
dimension.  

More precisely, let $\Delta(K)$ denote the algebraically shifted
complex of simplicial complex $K$, and let $\rbeti{j}(K,L)=\dimk
\rhomi{j}(K,L;\kk)$ be the dimension of the $j$th reduced relative
homology group over a field $\kk$ of a pair of simplicial complexes $L
\subseteq K$.  Then
$\rbeti{j}(K,L) \leq \rbeti{j}(\Delta(K),\Delta(L))$ for all $j$.

The theorem is motivated by somewhat similar results about Gr\"obner
bases and generic initial ideals.  Parts of the proof use Gr\"obner
basis techniques.
\end{abstract}

\maketitle

\section{Introduction}\label{se:intro}
Algebraic shifting is a remarkable procedure that finds, for any
simplicial complex $K$, a shifted (and hence combinatorially simpler)
simplicial complex $\Delta(K)$ with many of the same properties as
$K$.  For instance, the $f$-vector and homology Betti numbers are
preserved; Bj\"orner and Kalai~\cite{BjKal:1} used this fact to
characterize the $f$-vectors and Betti numbers of simplicial
complexes.

However, the situation for pairs of complexes and relative homology is
different.  In a simple example on three vertices
(Example~\ref{ex:pardue}), algebraically shifting a pair of complexes
increases their relative homology in dimensions $0$ and $1$.  Upon
seeing this one example, Keith Pardue (private communication)
conjectured that algebraic shifting always weakly increases relative
homology in every dimension.  Our main result (Theorem~\ref{th:main})
is that this conjecture is true.

Pardue's conjecture was grounded in more than just this one simple
example.  Algebraic shifting, which takes place in exterior
(anti-commutative) algebra, is similar to using Gr\"obner bases and
generic initial ideals in commutative algebra (see
Section~\ref{se:alg.shift}).  Quantities such as free resolution Betti
numbers weakly increase upon taking generic initial ideals (see,
\eg,~\cite{Hulett:thesis, Hulett, Bigatti}, and
Section~\ref{se:alg.shift}).  Pardue's insight was that these results
would carry over to algebraic shifting.

It would be ideal, then, to prove his conjecture by translating the
algebraic shifting problem to a generic initial ideal problem, and
then invoking the existing results.  However, this approach has been
unsuccessful, so far.  The proof here, while motivated at points by
Gr\"obner basis ideas (see Lemma~\ref{th:gamma.exist}), instead
directly refines Bj\"orner and Kalai's correspondence between the
homology of the original complexes and the combinatorics of the
algebraically shifted complexes.  The hope is that this result will
serve as further evidence of the deeper
connection between algebraic shifting and generic initial ideals.

Background and notation on simplicial complexes, including homology,
shifted complexes, and near-cones is in Section~\ref{se:scx}.
Algebraic shifting is reviewed and compared to generic initial ideals
in Section~\ref{se:alg.shift}.  In Section~\ref{se:rel.hom}, we use
Gr\"obner basis ideas to define a nice basis of a space associated
with a pair of complexes $(K,L)$, and then use this basis to compare
key components of the homology groups of $(K,L)$ and
$(\Delta(K),\Delta(L))$.  We prove our main result
(Theorem~\ref{th:main}) in Section~\ref{se:proof}.

\section{Simplicial complexes}\label{se:scx}
For basic definitions of simplicial complexes and their homology and
relative homology, see, \eg,~\cite[Chapter~1]{Mu}
or~\cite[Section~0.3]{St:CCA}.  We allow the empty simplicial complex
$\emptyset$ consisting of no faces; all other complexes must contain
the empty set as a $(-1)$-dimensional face.  We also allow the complex
$\{\emptyset\}$ consisting of only the empty face, but we do
distinguish between the two complexes $\emptyset$ and $\{\emptyset\}$.
Let $K_j$ denote the set of $j$-dimensional faces of a simplicial
complex $K$.  The \dfn{$f$-vector} of $K$ is the sequence
$(f_0,\ldots,f_{d-1})$, where $f_j = f_j(K) = \abs{K_j}$ and $d-1 =
\dim K$.  The same notion of $f$-vector will apply in this paper to
every finite collection of sets.

Let $\kk$ be a field, fixed throughout the paper.  The $j$th
\dfn{Betti number} of a simplicial complex $K$ is
$\rbeti{j}=\rbeti{j}(K)=\dimk\rhomi{j}(K)$, where $\rhomi{j}(K)$ is
the $j$th reduced homology group of $K$ (with respect to $\kk$).
Similarly, the $j$th \dfn{relative Betti number} of a pair of
simplicial complexes $L
\subseteq K$ is $\rbeti{j}=\rbeti{j}(K,L)=\dimk\rhomi{j}(K,L)$, where
$\rhomi{j}(K,L)$ is the $j$th reduced \dfn{relative homology group} of
the pair $(K,L)$ (with respect to $\kk$).

``Reduced'' homology means precisely to treat the empty set as a face
of any non-empty complex, so $\rbeti{0}$ is one less than the number
of connected components of $\Delta$, and hence one less than the
``unreduced'' $\rbeti{0}$.  Furthermore, $\rbeti{-1}=0$, unless
$\Delta=\{\emptyset\}$, in which case $\rbeti{-1}=1$.
Reduced relative homology, which also treats the empty set as a face
of any non-empty complex, is the same as unreduced
relative homology, except that $\rbeti{-1}(\{\emptyset\},\emptyset)=1$;
for any other pair of complexes,
$\rbeti{-1}=0$.

%Recall that over a field $\kk$, $\dimk \rcohomi{j}(K;\kk) = 
%\dimk \rhomi{j}(K;\kk)$, so that Betti numbers measure reduced homology
%as well as reduced cohomology.

\begin{defn}
If $S=\setseq{s}{1}{j}$ and $T=\setseq{t}{1}{j}$ are $j$-subsets
of integers,
then:
\begin{itemize}
\item $S \leq_P T$ under the standard \dfn{partial order} if $s_p
  \leq t_p$ for all $p$; and
\item $S <_L T$ under the \dfn{lexicographic order} if there is a
  $q$ such that $s_q < t_q$ and $s_p = t_p$ for $p < q$ .
\end{itemize}
\end{defn}

Lexicographic order is a total order which  refines the 
partial order.

\begin{defn}
A collection $\collC$ of $k$-subsets is \dfn{shifted} if $S \leq_P T$
and $T \in \collC$ together imply that $S \in \collC$.  A simplicial
complex 
$\Delta$ is \dfn{shifted} if the set of $j$-dimensional faces of
$\Delta$ is shifted for every $j$.
\end{defn}

Bj\"orner and Kalai showed in~\cite{BjKal:1} that
shifted complexes are near-cones, which we now define.

\begin{defn}
A \dfn{near-cone} with \dfn{apex} $v_0$ is a simplicial complex
$\Delta$ satisfying the following property:
For all $F \in \Delta$, if $v_0 \not\in F$ and $w \in F$, then 
$$(F - \{w\}) \cup \{v_0\} \in \Delta.$$  
For $\Delta$ a near-cone with apex $v_0$, let
$B(\Delta) = 
	\{F \in \Delta\colon F \cup \{v_0\} \not\in \Delta\}$
and 
$\Delta' 
  = \{F \in \Delta\colon v_0 \not\in F, F \cup \{v_0\} \in \Delta\}$;
then 
$$\Delta = (v_0 \ast \Delta') \disun B(\Delta),$$
where $\ast$ denotes topological join (so
$v_0 \ast \Delta' = 
  \Delta' \disun \{\{v_0\} \disun F\colon F \in \Delta'\}$).
Both $\Delta'$ and $\Delta' \disun B(\Delta)$ are subcomplexes of $\Delta$.
If $B(\Delta) = \emptyset$, then $\Delta$ is simply a \dfn{cone}.
\end{defn}

Note, in particular, that $\emptyset$ and $\{\emptyset\}$ are
near-cones (the condition in the definition is vacuous in this case)
and that $\emptyset = v_0 \ast \emptyset$ and 
$\{\emptyset\} = (v_0 \ast \emptyset) \disun \{\emptyset\}$.  
If $\Delta$ is a near-cone with apex $v_0$, then $v_0$ is one of the
vertices of $\Delta$, unless $\Delta=\emptyset\ {\rm or}\ \{\emptyset\}$.

It is not hard to see that shifted simplicial complexes are near-cones
with apex $1$.

Every $F \in B(\Delta)$ is maximal in $\Delta$, so the collection of
faces in $B(\Delta)$ forms an antichain.  Further, $f(B(\Delta)) =
\betat(\Delta)$, which follows by contracting $v_0 \ast \Delta'$ to
$v_0$, leaving a sphere for every face in
$B(\Delta)$~\cite[Theorem~4.3]{BjKal:1}.
In other words, 
if $\Delta$ is a near-cone with apex $v_0$, then
\begin{equation}\label{eq:beta.shift}
\rbeti{j}(\Delta)= \abs{\{F\in \Delta_j\colon 
	v_0 \not\in F, \nd v_0 \disun F \not\in \Delta\}}.
\end{equation}
This observation is generalized by~\cite[Lemma~8]{art:relative}:
If $\Gamma \subseteq \Delta$ is a pair of near-cones with common apex
$v_0$, then
\begin{multline}\label{eq:pair.near-cones}
\rbeti{j}(\Delta,\Gamma)=
  \abs{\{F\in (\Delta - \Gamma)_j\colon 
	v_0 \not\in F, \nd v_0 \disun F \not\in \Delta\}}\\  +
  \abs{\{G\in (\Delta - \Gamma)_j\colon 
	v_0 \in G, \nd G - \{v_0\} \in \Gamma\}}.
\end{multline}

In light of the formulation of the homology of near-cones that
equation~\eqref{eq:beta.shift} gives,
equation~\eqref{eq:pair.near-cones} is approximately the near-cone
equivalent of using the long exact sequence
(e.g.,~\cite[Theorem~23.3]{Mu})
\begin{equation}\label{eq:les}
\dots \rightarrow \rhomi{j}(L) \xrightarrow{i_*} \rhomi{j}(K) 
      \xrightarrow{\pi_*} \rhomi{j}(K,L) 
      \xrightarrow{\partial_*} \rhomi{j-1}(L) \rightarrow \cdots
\end{equation}
to compute
$$\rbeti{j}(K,L) = \dim(\im (\pi_*)_j) + \dim(\im (\partial_*)_j)$$
for an arbitrary pair $(K,L)$.

\section{Algebraic shifting}\label{se:alg.shift}

Algebraic shifting transforms a simplicial complex into a shifted
simplicial complex with the same $f$-vector and Betti numbers.  
It also preserves many
algebraic properties of the original complex.  Algebraic shifting was
introduced by Kalai in~\cite{Kal:Eck}; our exposition is summarized
from~\cite{BjKal:1} and included for completeness (see
also~\cite{BjKal:NYAS,Kal:AS}).  We start with the exterior face ring.

\begin{defn} 
Let $K$ be a $(d-1)$-dimensional 
simplicial complex with vertices $V=\set{e}{1}{n}$ 
linearly ordered $\ltseq{e}{1}{n}$.  Let $\Lambda(\kV)$ denote the 
exterior algebra of the vector space $\kV$; it has a $\kk$-vector 
space basis consisting of all the monomials 
$e_S := \wedgesubseq{e}{i}{1}{j}$, 
where $S=\setsubseq{e}{i}{1}{j} \subseteq V$ 
(and $e_{\emptyset}=1$).  
Note that $\Lambda(\kV)=\oplus_{j=0}^n \Lambda^j(\kV)$ is a graded
$\kk$-algebra, and that $\Lambda^j(\kV)$ has basis 
$\{e_S\colon \abs{S}=j\}$.
Let $(I_K)_j$ be the subspace of $\Lambda^{j+1}(\kV)$ generated by the
basis $\{e_S\colon \abs{S}=j+1, \nd  S \not\in K\}$.
Then $I_K:=\oplus_{j=-1}^{d-1} (I_K)_j$ is the homogeneous
graded ideal 
of $\Lambda(\kV)$ generated by $\{e_S\colon S \not\in K\}$.  
Let $\Lambda_j[K] := \Lambda^{j+1}(\kV)/(I_K)_j$.
Then the graded quotient algebra 
$\Lambda[K]:= \oplus_{j=-1}^{d-1} \Lambda_j[K] 
	         = \Lambda(\kV)/I_K$
is called the \dfn{exterior face ring} of $K$ (over $\kk$).  
\end{defn} 

The exterior face ring is the exterior algebra analogue to the
Stanley-Reisner face ring of a simplicial complex~\cite{St:CCA}.  For
$x\in\kV$, let $\imt{x}$ denote the image of $x$ in $\Lambda[K]$.

\begin{defncite}{Kalai}
  Let $\set{f}{1}{n}$ be a ``generic'' basis of $\kV$, \ie,
  $$f_i=\sum^n_{j=1}\alpha_{ij}e_j,$$ where the $\alpha_{ij}$'s are
  $n^2$ transcendentals, algebraically independent over $\kk$.
  Define $f_S:=\wedgesubseq{f}{i}{1}{j}$ for
  $S=\setseq{i}{1}{j}$ (and set $f_\emptyset = 1$).  Let
$$\Delta(K,\kk) :=
        \{S \subseteq [n]\colon \imf{S} \not\in
                \spn \{\imf{R}\colon R <_L S \}\}$$ 
be the \dfn{algebraically shifted complex} obtained from $K$;
we will write $\Delta(K)$ instead of $\Delta(K,\kk)$ when the field is
understood to be $\kk$.  In other words, the $j$-subsets of
$\Delta(K)$ can be chosen by listing all the $j$-subsets of $[n]$ in
lexicographic order and omitting those that are in the span of earlier
subsets on the list, modulo $I_K$ and with respect to the $f$-basis.

The algebraically shifted complex $\Delta(K)$ is (as its name suggests)
shifted, and is independent of the numbering of the vertices of 
$K$~\cite[Theorem~3.1]{BjKal:1}.
\end{defncite} 

It is easy to see that algebraic shifting preserves the $f$-vector,
\ie, $f_j(K)=f_j(\Delta(K))$.  Bj\"orner and Kalai~\cite{BjKal:1} showed that
algebraic shifting also preserves Betti numbers, 
\ie, $\rbeti{j}(K)=\rbeti{j}(\Delta(K))$.  The reason lies in the relation between
algebraic shifting and coboundaries.  Define the \dfn{weighted coboundary operator}
$\delta\colon \Lambda[K]\maps\Lambda[K]$ by
$\delta(x)=\imf{1} \wedge x$, so
$$\delta(\ime{S})= \imf{1} \wedge \ime{S} =
\sum^n_{j=1}\alpha_{1j} \ime{j} \wedge \ime{S} =
\sum_{\condns{j\not\in S}{S \cup \{j\} \in K}}
        \pm \alpha_{1j} \ime{S \cup \{j\}}$$ 
(hence the name weighted coboundary operator).
Betti numbers may be computed using
this $\delta$, \ie,
$\rbeti{j}(K)=\dimk(\ker \delta)_j/(\im \delta)_j$~\cite[pp.~289--290]{BjKal:1}.  
Furthermore, the action of $\delta$ 
on many members of the $\imt{f}$-basis
is easy to describe: $\delta(\imf{F})$ equals $\imf{1\disun F}$ if $1
\not\in F$ and $1 \disun F \in \Delta(K)$, but is zero if $1\in F$
(the third case, when $1 \not\in F$, but $1 \disun F \not\in
\Delta(K)$, is harder, and we shall not need it).

What about relative homology?  First, we note a result of
Kalai's~\cite[Theorem~2.2]{Kal:AS} 
that if $L \subseteq K$ is a pair of simplicial
complexes, then $\Delta(L) \subseteq \Delta(K)$.  For every pair
$(K,L)$, we may then consider the pair $(\Delta(K),\Delta(L))$.  In
contrast to the single complex case, however, the homology of $(K,L)$
and $(\Delta(K),\Delta(L))$ need not coincide, as the following
example shows.

\begin{ex}\label{ex:pardue} 
Let $K$ be the simplicial complex on vertices $\{1,2,3\}$ whose
maximal faces are $\{1\}$ and $\{2,3\}$, and $L$ be the subcomplex
consisting of just the two vertices $\{1\}$ and $\{2\}$.  The only
shifted complex with three vertices and one edge has maximal faces
$\{1,2\}$ and $\{3\}$, so this must be $\Delta(K)$.  Furthermore, $L$
is the only simplicial complex with two vertices and no edges, so
$\Delta(L)=L$.  (See
Figure~\ref{fg:pardue}.)  
\begin{figure}
\begin{center}
\setlength{\unitlength}{0.012500in}%
\begin{picture}(458,36)(71,720)
\thicklines
\put(451,740){\circle{8}}
\put(525,740){\circle*{8}}
\put(364,740){\circle{8}}
\put(368,740){\line( 1, 0){ 79}}
\put(364,747){\makebox(0,0)[b]{$1$}}
\put(451,747){\makebox(0,0)[b]{$2$}}
\put(525,747){\makebox(0,0)[b]{$3$}}
\put(448,720){\makebox(0,0)[b]{$(\Delta(K),\Delta(L))$}}
\put(149,740){\circle{8}}
\put( 75,740){\circle{8}}
\put(236,740){\circle*{8}}
\put(232,740){\line(-1, 0){ 79}}
\put(149,747){\makebox(0,0)[b]{$2$}}
\put( 75,747){\makebox(0,0)[b]{$1$}}
\put(236,747){\makebox(0,0)[b]{$3$}}
\put(148,720){\makebox(0,0)[b]{$(K,L)$}}
\end{picture}
\end{center}
\caption{Example~\ref{ex:pardue}}\label{fg:pardue}
\end{figure}
But then it is easy to see that
$(\Delta(K),\Delta(L))$ has non-trivial relative homology in
dimensions $0$ and $1$, while $(K,L)$ has no non-trivial relative
homology.

Thus, the relative Betti numbers of $(\Delta(K),\Delta(L))$ are all at
least as large as those of $(K,L)$.  Theorem~\ref{th:main} shows that
this is true for any pair $(K,L)$.
\end{ex}

Algebraic shifting is the exterior algebra analogue of generic initial
ideals and Gr\"obner bases in commutative algebra, in the following
way.  If $I_K$ were instead a monomial ideal of a {\em polynomial}\
ring, then the algorithm used to create the list of non-faces of
$\Delta(K)$ would instead create a list of monomials generating the
generic initial ideal of $I_K$, denoted $\Gin(I_K)$.  For further
details of generic initial ideals, see, for
instance~\cite[Section~15.9]{Eisenbud}.  For more about the
relationship between generic initial ideals and algebraic shifting,
see~\cite{HeTe}.  For a more general exterior algebra version of
Gr\"obner bases and generic initial ideals,
see~\cite[Section~1]{AHH:Gotzmann}.

Theorem~\ref{th:main} bears some resemblance to results about generic
initial ideals (Section~\ref{se:alg.shift}).  For instance,
Hulett~\cite[proof of Lemma~1.24]{Hulett:thesis},~\cite[p.~2338]{Hulett} 
and Bigatti~\cite[proof of Theorem~3.7]{Bigatti} have shown that for
any homogeneous ideal $I$ in a polynomial ring, the free resolution
Betti numbers of its generic initial ideal $\Gin(I)$ are at least as
large as those of $I$.

\section{Relative homology}\label{se:rel.hom}

In order to say anything about $(\Delta(K),\Delta(L))$, we must first
consider $(K,L)$.  For $Q=K-L$ (the ``$Q$'' is for ``quotient''), we
define $\Delta(Q)=\Delta(K)-\Delta(L)$.  This is primarily a
combinatorial definition, with the algebra hidden in the computation
of $\Delta(K)$ and $\Delta(L)$.  We now examine how to interpret
$\Delta(Q)$ algebraically.  Let
$$\imt{Q}=\spn\{\ime{F}\colon F \in Q\}.$$
It is not hard to see, then, that we may algebraically shift
the subcomplex $L$ using $\Lambda[K]$ instead of $\Lambda[L]$, by
modding out by $\imt{Q}$ on $\Lambda[K]$ instead of by $I_L$ on
$\Lambda[L]$, since $\imt{Q}=\imt{I_L}$ (see~\cite[Section~3]{art:relative}).

Lemmas~\ref{th:gamma.exist} and~\ref{th:gamma.basis} show 
how $\Delta(Q)$ is related to $\imt{Q}$, namely
that $\Delta(Q)$ indexes a nice basis of $\imt{Q}$;
the construction is motivated by Gr\"obner basis ideas.  Then, guided
by earlier results about $\imt{Q}$ (summarized here as
Lemma~\ref{th:im.ker}), we use this basis of $\imt{Q}$ in
Lemmas~\ref{th:Sigma.top} and~\ref{th:Sigma.bottom} to compare key
subspaces of $\rhom(K,L)$ and $\rhom(\Delta(K),\Delta(L))$.

\begin{lemma}\label{th:gamma.exist}
If $F \in \Delta(Q)$, then there is a unique linear combination
$\sum_{\condns{G <_L F}{G \in \Delta(L)}} a_G \imf{G}$, such that 
$\imf{F}-\sum_{\condns{G <_L F}{G \in \Delta(L)}} a_G \imf{G} \in
\imt{Q}$.
\end{lemma}

\begin{proof}
Since $F \in \Delta(Q)$, and hence $F \not\in \Delta(L)$, we have
\begin{equation}\label{eq:normal.form}
\imf{F}-\sum_{G <_L F} a_G \imf{G} \in I_L = \imt{Q}
\end{equation}
for some $a_G$.  We may iterate this process on the $\imf{G}$'s for
which $G \not\in \Delta(L)$, replacing {\em them} by lexicographically
earlier linear combinations that are equal modulo $I_L=\imt{Q}$
until every $G$ in equation~\eqref{eq:normal.form} is in 
$\Delta(L)$.  
This eventually terminates, since lexicographic order is a total order.
(In Gr\"obner basis theory, this procedure is known as finding the
normal form~\cite[Definition~2.1.3]{AdLou}.)

To show these $a_G$ are unique, assume that also
\begin{equation}\label{eq:gamma.unique.2}
\imf{F}-\sum_{\condns{G <_L F}{G \in \Delta(L)}} b_G \imf{G} \in
\imt{Q}.
\end{equation}
Then by subtracting equation~\eqref{eq:gamma.unique.2} from
equation~\eqref{eq:normal.form}, we get
\begin{equation}\label{eq:gamma.unique.3}
\sum_{\condns{G <_L F}{G \in \Delta(L)}} (b_G - a_G) \imf{G} \in
\imt{Q}.
\end{equation}

If any $b_G - a_G$ in equation~\eqref{eq:gamma.unique.3} is non-zero,
let $G_0$ index the lexicographically last of these; then 
$$(a_{G_0}-b_{G_0})\imf{G_0}
 -\sum_{\condns{G <_L G_0}{G \in \Delta(L)}} (b_G-a_G) \imf{G} \in
\imt{Q},$$ and $$\imf{G_0} -\sum_{\condns{G <_L G_0}{G \in
\Delta(L)}} \left(\frac{b_G-a_G}{a_{G_0}-b_{G_0}}\right) \imf{G} \in
\imt{Q},$$
which contradicts $G_0 \in \Delta(L)$.
\end{proof}

\begin{defn}
By Lemma~\ref{th:gamma.exist}, we may define, for any $F \in
\Delta(Q)$, $\img{F}$ supported on $\Delta(L)$ such that 
$\imf{F}-\img{F} \in \imt{Q}$.
\end{defn}

\begin{lemma}\label{th:gamma.basis}
$\{\imf{F}-\img{F}\colon F \in \Delta(Q)\}$ is a basis of $\imt{Q}$.
\end{lemma}

\begin{proof}
We first show that $\{\imf{F}-\img{F}\colon F \in \Delta(Q)\}$ is
linearly independent.  Assume otherwise; 
\begin{equation}\label{eq:gamma.ind}
\sum_{F \in \Delta(Q)} b_F(\imf{F}-\img{F})=0,
\end{equation}
where $b_{F_0} \neq 0$ for some $F_0 \in \Delta(Q)$.  When expanding
the sum on the left-hand side of equation~\eqref{eq:gamma.ind} in the
$\{\imf{F}\colon F \in \Delta(K)\}$ basis, the coefficient of
$\imf{F_0}$ will be $b_{F_0} \neq 0$, since the $\img{F}$ are all
supported on $\Delta(L)$, and so cannot cancel $\imf{F_0}$.  So
$\{\imf{F}-\img{F}\colon F \in \Delta(Q)\}$ is a set of $\abs{\Delta(Q)}$
linearly independent vectors in $\imt{Q}$.

On the other hand, $\imt{Q}$ is a
$\abs{K-L}=\abs{\Delta(Q)}$-dimensional vector space, so
$\{\imf{F}-\img{F}\colon F \in \Delta(Q)\}$ must be a basis.
\end{proof}

Now we see how $\imt{Q}$ can help compute homology and relative homology.
We adopt the shorthand $$\delta^{-1}\imt{Q}=
\{\imt{x} \in \Lambda[K]\colon \delta \imt{x} \in \imt{Q}\}.$$

\begin{lemma}\label{th:im.ker}
For any pair of simplicial complexes $L \subseteq K$,
\begin{alph-list}
\item\label{it:Gamma}
$\rbeti{j}(L)=\dim 
((\delta^{-1}\imt{Q})/(\im \delta + \imt{Q}))_j$; and
\item\label{it:Sigma} 
$\rbeti{j}(K,L)=\dim
((\ker \delta \cap \imt{Q})/(\delta \imt{Q}))_j$.
\end{alph-list}
\end{lemma}

\begin{proof}
This is~\cite[Lemmas 2 and 4]{art:relative}, where the notation
$\Lambda[\Sigma]$ was used in place of $\imt{Q}$.
\end{proof}

Lemma~\ref{th:im.ker}\ref{it:Sigma} suggests that in order to compute
$\rbeti{j}(K,L)$, we examine $\ker \delta \cap \imt{Q}$ and
$\delta\imt{Q}$.  However, $\im \delta \cap \imt{Q}$ turns out to be
easier to handle than $\ker \delta \cap \imt{Q}$.  The next two lemmas
compare $\im \delta \cap \imt{Q}$ and $\delta \imt{Q}$ to
subspaces of $\imt{Q}$ indexed by combinatorially defined sets of
$\Delta(Q)$. These two comparisons will combine to prove the key
inequality in the proof of Theorem~\ref{th:main}.

\begin{lemma}\label{th:Sigma.top}
$\im \delta \cap \imt{Q} \subseteq 
	\spn \{\imf{1 \disun F}-\img{1 \disun F}\colon
		1 \not\in F,\nd 1 \disun F \in \Delta(Q)\}$.
\end{lemma}

\begin{proof}
Let $\imt{x} \in \im \delta \cap \imt{Q}$.  By
Lemma~\ref{th:gamma.basis}, we can write
\begin{equation}\label{eq:supp.Sigma}
\imt{x} = \sum_{G \in \Delta(Q)} a_G (\imf{G}-\img{G})
\end{equation}
uniquely, since $\imt{x}\in\imt{Q}$.  Similarly, 
by~\cite[equation~(3.5)]{BjKal:1}, we can also write
\begin{equation}\label{eq:supp.1}
%\imt{x}=\sum_{1 \disun F \in \Delta(K)} b_F \imf{1 \disun F}
\imt{x}=\sum_{\condns{1 \not\in F}{1 \disun F \in \Delta(K)}} 
		b_F \imf{1 \disun F}
\end{equation}
uniquely, since $\imt{x}\in\im \delta$.  Now, by definition, the
support of $\img{G}$ is entirely on $\Delta(L)$.  Of course, the
support of $\imt{x}$ in $\Delta(Q)$ must be the same in
equations~\eqref{eq:supp.Sigma} and~\eqref{eq:supp.1},
so we must be able to write every $G \in \Delta(Q)$ such that $a_G
\neq 0$ as $G=1 \disun F$ for some $F$.  Therefore,
equation~\eqref{eq:supp.Sigma} can be rewritten as
$$\imt{x} = \sum_{\condns{1 \not\in F}{1 \disun F \in \Delta(Q)}} 
		a_{1 \disun F} (\imf{1 \disun F}-\img{1 \disun F}),$$
implying the lemma.
\end{proof}

\begin{lemma}\label{th:Sigma.bottom}
There is a subspace $\imt{Q'}$ of $\imt{Q}$ such that
$$\dim \delta \imt{Q'} = \abs{ \{F \in \Delta(Q)\colon
	1 \not\in F, \nd 1 \disun F \in \Delta(Q)\}}.$$
\end{lemma}

\begin{proof}
Let $\imt{Q'}=\spn \{\imf{F}-\img{F}\colon
		F \in \Delta(Q), \nd 1 \not\in F, \nd
		1 \disun F \in \Delta(Q)\}$.  
First note that by definition of $\img{F}$ (and $F \in \Delta(Q)$),
each $\imf{F}-\img{F} \in \imt{Q}$, so $\imt{Q'}$ is a subspace of
$\imt{Q}$.  Clearly, we only need to show that 
$\{\delta(\imf{F}-\img{F})\colon
		F \in \Delta(Q), \nd 1 \not\in F, \nd
		1 \disun F \in \Delta(Q)\}$
is linearly independent.

By definition of $\img{F}$ and Lemma~\ref{th:gamma.exist}, we may
write each
$$\imf{F}-\img{F}=\imf{F}-
  \sum_{\condns{G<_L F}{G\in\Delta(L)}} b_{F,G}\imf{G}$$
for some $b_{F,G}$'s.
Furthermore, we are assuming $1 \disun F \in \Delta(Q)$ for each $F$,
so $1 \disun F \in \Delta(K)$, and thus $\delta \imf{F}=\imf{1 \disun F}$.
For each $G <_L F$, if $1 \in G$, then $\delta\imf{G}=0$;  otherwise 
$1 \disun G <_L 1 \disun F \in \Delta(K)$, so $1 \disun G \in
\Delta(K)$, and $\delta \imf{G}=\imf{1\disun G}$.  Therefore
$$\delta(\imf{F}-\img{F})=\imf{1 \disun F}-
  \sum_{\condns{G<_L F}{1 \not\in G, \nd G\in\Delta(L)}} 
		b_{F,G}\imf{1 \disun G}.$$

To show that $\{\delta(\imf{F}-\img{F})\colon
		F \in \Delta(Q), \nd 1 \not\in F, \nd
		1 \disun F \in \Delta(Q)\}$
is linearly independent, assume
\begin{equation}\label{eq:Sigma.ind}
0=\sum_{\condns{F \in \Delta(Q)}{1 \not\in F,\nd 1 \disun F \in
\Delta(Q)}} c_F(\imf{1 \disun F} - 
  \sum_{\condns{G <_L F}{1 \not\in G, \nd G \in \Delta(L)}}
    b_{F,G}\imf{1 \disun G}).
\end{equation}
Now, $f_{1 \disun G}$'s appearing in equation~\eqref{eq:Sigma.ind} all
satisfy $G\in \Delta(L)$, while all the $\imf{1 \disun F}$'s satisfy 
$F \in \Delta(Q)$, so there is no cancellation between the 
$\imf{1 \disun G}$'s and the $\imf{1 \disun F}$'s.  But all the 
$\imf{1 \disun F}$'s are distinct members of the 
$\{\imf{H}\colon H \in \Delta(K)\}$ basis, so there is no cancellation
among the $\imf{1 \disun F}$'s.  Therefore all the $c_F$'s must be
zero, and the $\delta(\imf{F}-\img{F})$'s are linearly independent.
\end{proof}

\section{Proof of main theorem}\label{se:proof}
We start with an easy lemma.
\begin{lemma}\label{th:KIJ}
If $I$, $J$, and $K$ are subspaces of a vector space and $I \subseteq
K$, then
$$\dim (K/I) 
= \dim ((K \cap J)/(I \cap J))
+ \dim ((K + J)/(I + J)).$$
\end{lemma}

\begin{proof}
Simply expand the right-hand side as
\begin{multline*}
(\dim (K \cap J)-\dim (I \cap J))
+((\dim K + \dim J - \dim (K \cap J)) \\ - 
  (\dim I + \dim J - \dim (I \cap J)))
\end{multline*}
by the standard vector space argument
$\dim (A + B) = \dim A + \dim B - \dim (A \cap B)$, applied twice.
This expression then easily simplifies to 
$\dim K - \dim I = \dim (K/I)$.
\end{proof}

\begin{thm}\label{th:main}
For any pair of simplicial complexes $L \subseteq K$,
$$\rbeti{j}(K,L) \leq \rbeti{j}(\Delta(K),\Delta(L))$$
for all $j$.
\end{thm}

\begin{proof} 
Because $\Delta(K)$ and $\Delta(L)$ are shifted, and therefore
near-cones with apex $1$, we may use equations~\eqref{eq:beta.shift}
and~\eqref{eq:pair.near-cones} to compute the homology of $\Delta(K)$,
$\Delta(L)$, and  $(\Delta(K),\Delta(L))$.
The sets in these equations overlap in a nice way.  In particular, if
we let
\begin{align*}
C_{KQ} &= \{F \in \Delta(Q)\colon 1 \not\in F,\nd 
		1 \disun F \not\in \Delta(K)\},\\
C'_{LQ} &= \{G \in \Delta(Q)\colon 1 \in G,\nd 
		G-1 \in \Delta(L)\},\\
C_{LQ} &= \{F \in \Delta(L)\colon 1 \not\in F,\nd 
		1 \disun F \in \Delta(Q)\},\ {\rm and}\\
%\intertext{and}
C_{KL} &= \{F \in \Delta(L)\colon 1 \not\in F,\nd 
		1 \disun F \not\in \Delta(K)\},
\end{align*}
then it is not hard to see, from equation~\eqref{eq:beta.shift}, that
\begin{align*}
\rbeti{j}(\Delta(K))&=\abs{(C_{KQ})_j}+\abs{(C_{KL})_j},\\
\rbeti{j}(\Delta(L))&=\abs{(C_{LQ})_j}+\abs{(C_{KL})_j},\\
\intertext{and, from equation~\eqref{eq:pair.near-cones}, that}
\rbeti{j}(\Delta(K),\Delta(L))&=\abs{(C_{KQ})_j}+\abs{(C'_{LQ})_j}.\\
\end{align*}
(We name these sets ``$C$'' because they are combinatorial.)
An easy bijection ($F \leftrightarrow 1 \disun F = G$) shows that
$$\abs{(C_{LQ})_j}=\abs{(C'_{LQ})_{j+1}}.$$
Continuing the analogy begun at the end of Section~\ref{se:scx}
between formulas for homology of near-cones and the long exact
sequence~\eqref{eq:les}, $C_{KQ}$ corresponds to $\im \pi_*$, 
$C_{LQ}$ and $C'_{LQ}$ correspond to $\im \partial_*$, and
$C_{KL}$ corresponds to $\im i_*$.

We can find ``corresponding'' subspaces in $\Lambda[K]$; define
\begin{align*}
A_{KQ}&=(\ker \delta \cap \imt{Q})/(\im \delta \cap \imt{Q}),\\
A'_{LQ}&=(\im \delta \cap \imt{Q})/(\delta \imt{Q}),\\
A_{LQ}&=(\delta^{-1}\imt{Q})/(\ker \delta + \imt{Q}),\ {\rm and}\\
%\intertext{and}
A_{KL}&=(\ker \delta + \imt{Q})/(\im \delta + \imt{Q}).
\end{align*}
(We name these spaces ``$A$'' because they are algebraic.)
Then by Lemma~\ref{th:KIJ}, 
\begin{align*}
\rbeti{j}(K)&=\dim (A_{KQ})_j + \dim (A_{KL})_j;\\
\intertext{by Lemma~\ref{th:im.ker},}
\rbeti{j}(L)&=\dim (A_{LQ})_j + \dim (A_{KL})_j,\ {\rm and}\\
\rbeti{j}(K,L)&=\dim (A_{KQ})_j +\dim (A'_{LQ})_j;
\end{align*}
and, by~\cite[Lemma~5]{art:relative},
$$ (A_{LQ})_j \cong (A'_{LQ})_{j+1}.$$

We will show how the dimension of each $A$ subspace compares with the
cardinality of the corresponding $C$ set with the same subscript.
Because algebraic shifting preserves homology,
\begin{align}
 \dim (A_{LQ})_j +\dim (A_{KL})_j = \rbeti{j}(L) &= \rbeti{j}(\Delta(L))
=\abs{(C_{LQ})_j}+\abs{(C_{KL})_j} \label{eq:new.1}\\
\intertext{and}
 \dim (A_{KQ})_j +\dim (A_{KL})_j = \rbeti{j}(K) &= \rbeti{j}(\Delta(K))
=\abs{(C_{KQ})_j}+\abs{(C_{KL})_j}. \label{eq:new.2}
\end{align}
By Lemma~\ref{th:Sigma.top}, 
\begin{align*}
\dim (\im \delta \cap \imt{Q})_{j+1} &\leq
  \abs{\{F \in \Delta(K)_{j}\colon 1 \not\in F, \nd 1 \disun F \in\Delta(Q) \}},\\
\intertext{and by Lemma~\ref{th:Sigma.bottom},}
\dim (\delta \imt{Q})_{j+1} &\geq \dim (\delta \imt{Q'})_{j+1}\\ &=
 \abs{\{F \in \Delta(Q)_{j}\colon 1 \not\in F, \nd 1 \disun F \in\Delta(Q) \}}.
\end{align*}
(The index shift, of $j+1$ to $j$, in the above inequalities arises
because the $(j+1)$-dimensional basis elements for $\im \delta \cap
\imt{Q}$ and $\delta\imt{Q'}$ are the $(j+1)$-dimensional elements
$\imf{1\disun F} - \img{1\disun F}$ and $\delta(\imf{F} - \img{F})$,
respectively, each of which is $(j+1)$-dimensional precisely when $F$
is $j$-dimensional.)
Since $\Delta(L)$ is the complement of $\Delta(Q)$ with respect to
$\Delta(K)$, then,
\begin{align*}%\label{eq:star}
\dim (A'_{LQ})_{j+1} &= \dim ((\im \delta \cap \imt{Q})/(\delta\imt{Q}))_{j+1}
= \dim (\im \delta \cap \imt{Q})_{j+1} - \dim (\delta\imt{Q})_{j+1}
					\notag\\
 &\leq \abs{\{F \in \Delta(L)_{j}\colon 
	1 \not\in F, \nd 1 \disun F \in\Delta(Q)\}} \notag\\
 &= \abs{(C_{LQ})_{j}},
\end{align*}
and so
$$\dim (A_{LQ})_j = \dim (A'_{LQ})_{j+1} \leq
                    \abs{(C_{LQ})_j} = \abs{(C'_{LQ})_{j+1}}$$

%\begin{align*}
%\dim (A_{LQ})_j &\leq \abs{(C_{LQ})_j}. \notag\\
%\intertext{Equation~\eqref{eq:new.1} then implies}
Equation~\eqref{eq:new.1} then implies
\begin{align*}
\dim (A_{KL})_j &\geq \abs{(C_{KL})_j},\\
\intertext{so by equation~\eqref{eq:new.2},}
\dim (A_{KQ})_j &\leq \abs{(C_{KQ})_j},\\
\intertext{and so finally}
\rbeti{j}(K,L)=\dim (A_{KQ})_j +\dim (A'_{LQ})_j &\leq 
               \abs{(C_{KQ})_j}+\abs{(C'_{LQ})_j} 
  = \rbeti{j}(\Delta(K),\Delta(L)).
\end{align*}
\end{proof}

\subsection*{Acknowledgements}
I gratefully acknowledge help from conversations with Keith Pardue,
Frank Sottile, and Eric Babson.

\end{document}